\documentclass[12pt]{article}

\usepackage[english]{babel}
\usepackage{amsmath}
\usepackage{amssymb,amsfonts,amsthm}
\usepackage{tikz,color,pgf,graphicx,url}
\usetikzlibrary{calc}
\usetikzlibrary{decorations.pathreplacing}
\usepackage{enumerate}     
\usepackage{bbm}
\usepackage{comment}
\usepackage[textsize=scriptsize,colorinlistoftodos]{todonotes}
\newcommand{\diff}{\mathrm{d}}

\usepackage[normalem]{ulem}

\usepackage{xcolor} 

\usepackage{hyperref}
\usepackage{cleveref}
\hypersetup{
	colorlinks   = true, 
	urlcolor     = black, 
	linkcolor    = black, 
	citecolor   = black 
}

\usepackage[
top=3cm,
bottom=3cm,
inner=2.5cm,      
outer=2.5cm,
footskip=40pt     
]{geometry}

\newtheorem{theorem}{Theorem}[section]

\newtheorem{lemma}[theorem]{Lemma}
\newtheorem{proposition}[theorem]{Proposition}
\newtheorem{conjecture}[theorem]{Conjecture}

\newtheorem{question}[theorem]{Question}

\usepackage[
backend=biber,
doi=false,
isbn=false,
url=false,
eprint=true,
sorting=nyt,
maxbibnames=999,
giveninits,
sortcites,]{biblatex}
\renewbibmacro{in:}{%
	\ifentrytype{article}{}{\printtext{\bibstring{in}\intitlepunct}}}
\AtEveryBibitem{%
	\clearfield{number}%
	\clearfield{note}%
	\clearfield{series}%
	\clearfield{pagetotal}}
\DeclareFieldFormat[article,inbook,incollection,inproceedings,patent,
thesis,unpublished]{title}{\mkbibemph{#1}\isdot}
\DeclareFieldFormat{journaltitle}{#1\isdot}
\DeclareFieldFormat{booktitle}{#1\isdot}
\DeclareFieldFormat[article]{pages}{#1}
\defcounter{biburlnumpenalty}{3000}
\newcommand{\eps}{\varepsilon}

\newcommand{\deep}{\mathcal{D}}
\newcommand{\midpath}{\mathcal{M}}

\newcommand{\cov}{\mathrm{cov}}
\newcommand{\EE}{\mathbb{E}}
\newcommand{\prob}{\mathbb{P}}
\newcommand{\N}{\mathbb{N}}
\newcommand{\R}{\mathbb{R}}
\newcommand{\depth}{\mathrm{depth}}
\newcommand{\inv}{^{-1}}
\DeclareMathOperator{\Var}{Var}

\newcommand{\indic}[1]{\mathbbm{1}_{#1}}

\begin{document}

\title{Packing and finding paths in sparse random graphs}


\author{Vesna Ir\v si\v c \footnote{\href{mailto:vesna.irsic@fmf.uni-lj.si}{vesna.irsic@fmf.uni-lj.si}, Faculty of Mathematics and Physics, University of Ljubljana, Slovenia}
\and Julien Portier \footnote{\href{mailto:jp899@cam.ac.uk}{jp899@cam.ac.uk}, Department of Pure Mathematics and Mathematical Statistics (DPMMS), University of Cambridge, Wilberforce Road, Cambridge, CB3 0WA, United Kingdom}
\and Leo Versteegen \footnote{\href{mailto:lversteegen.math@gmail.com}{lversteegen.math@gmail.com}, Department of Pure Mathematics and Mathematical Statistics (DPMMS), University of Cambridge, Wilberforce Road, Cambridge, CB3 0WA, United Kingdom}}

\maketitle

\begin{abstract}
Let $G\sim G(n,p)$ be a (hidden) Erd\H{o}s-Rényi random graph with $p=(1+\eps)/n$ for some fixed constant $\eps >0$. Ferber, Krivelevich, Sudakov, and Vieira showed that to reveal a path of length $\ell=\Omega(\frac{\log(1/\eps)}{\eps})$ in $G$ with high probability, one must query the adjacency of $\Omega(\frac{\ell}{p\eps\log(1/\eps)})$ pairs of vertices in $G$, where each query may depend on the outcome of all previous queries. Their result is tight up to the factor of $\log(1/\eps)$ in both $\ell$ and the number of queries, and they conjectured that this factor could be removed. We confirm their conjecture. The main ingredient in our proof is a result about path-packings in random labelled trees of independent interest. Using this, we also give a partial answer to a related question of Ferber, Krivelevich, Sudakov, and Vieira. Namely, we show that when $\ell=o((t/\log t)^{1/3})$, the maximum number of vertices covered by edge-disjoint paths of length at least $\ell$ in a random labelled tree of size $t$ is $\Theta(t/\ell)$ with high probability.
\end{abstract}


\noindent
{\bf AMS Subj.\ Class.\ (2020):} 05C80, 05C85

\section{Introduction}
\label{sec:intro}
In the \emph{Erdős–Rényi random graph model} a graph $G \sim G(n,p)$ is obtained by independently adding an edge between each pair of $n$ labelled vertices with probability $p \in [0,1]$. We say that an event whose probability depends on $n$ occurs \emph{with high probability} (or \emph{whp} for short) if its probability goes to $1$ as $n$ goes to infinity.

In the companion papers \cite{ferber2016finding, ferber-2017} Ferber, Krivelevich, Sudakov, and Vieira introduced the following type of algorithmic problem on Erdős–Rényi random graphs.
Let $\mathcal{P}$ be an increasing graph property, i.e., a graph property such that when $G_1$ is a subgraph of $G_2$ and $G_1$ has the property $\mathcal{P}$, then $G_2$ does as well. Suppose that $p=p(n)$ is chosen such that when $G$ is sampled according to the Erd\H{o}s-Rényi graph model $G(n,p)$, then $G$ has the property $\mathcal{P}$ with high probability. 

An adaptive algorithm queries for pairs of vertices in $V(G)$ whether or not they are an edge of $G$. Each query is allowed to depend on the outcome of previous queries, and the algorithm terminates once the revealed edges witness that $G$ has $\mathcal{P}$. The principal objective in this setting is to determine the minimum number of queries that an adaptive algorithm requires to terminate with high probability. After the research of adaptive algorithms on random graphs was initiated in  \cite{ferber2016finding, ferber-2017}, problems of this type were studied in multiple papers, see for instance \cite{alweiss2021subgraph,feige2020finding,huleihel2024random, rashtchian2021average}.

If all graphs with property $\mathcal{P}$ have at least $m$ edges, then it is obvious that one needs to query at least $(1+o(1)) \frac{m}{p}$ pairs of vertices, as the algorithm requires at least $m$ positive responses in order to satisfy the target property.
In \cite{krivelevich2013phase}, where a question of this type appeared implicitly, Krivelevich and Sudakov showed that this trivial lower bound is tight in the supercritical regime if $\mathcal{P}$ is the property of having a giant component.
More precisely, they showed that there is an adaptive algorithm which finds a connected component of size at least $\eps n/2$ in $G\sim G(n,p)$ in at most $\eps n^2/2$ queries with high probability when $p=\frac{1+\eps}{n}$.
Similarly, Ferber, Krivelevich, Sudakov, and Vieira proved in \cite{ferber2016finding} that this trivial lower bound is also tight for the property of having a Hamilton cycle.
They established that when $p\geq \ln n+\ln \ln n+\omega(1)$ there exists an adaptive algorithm that finds a Hamiltonian cycle in $G\sim G(n,p)$ in at most $(1+o(1))\frac{n}{p}$ queries with high probability.
The naive bound, however, is not always tight. In \cite{ferber-2017}, Ferber, Krivelevich, Sudakov, and Vieira showed that if the property of interest is to have a path of length $\ell$, then one may need substantially more than $\ell/p$ queries.

\begin{theorem}[{\cite[Theorem 1]{ferber-2017}}]
\label{thm:FKSV}
 There exists an absolute constant $C >0$ such that the following holds. For every
constant $q \in (0,1)$, there exists $n_0,\eps_0$ such that for every fixed $\eps \in (0, \eps_0)$ and any $n \geq n_0$, there is no adaptive algorithm which reveals a path of length $\ell \geq \frac{3C}{\eps} \log(\frac{1}{\eps})$ with probability at least $q$ in $G \sim G(n,p)$, where $p=\frac{1+\eps}{n}$, by quering at most $\frac{q\ell}{8640Cp\eps \log(\frac{1}{\eps})}$ pairs of vertices.
\end{theorem}

They conjectured that the $\log(1/\eps)$ factor in the denominator could be removed, and noticed that this would be optimal, since an argument from \cite{krivelevich2013phase} can be adapted to show that there exists an adaptive algorithm which finds a path of length $\ell \leq \frac{1}{5}\eps^2n$ in $G(n,p)$ when $p=\frac{1+\eps}{n}$ with probability at least $1-\exp(-\Omega(\frac{\ell}{\eps}))$ in at most $O(\frac{\ell}{p\eps})$ queries.
In this paper, we confirm their conjecture in form of the following theorem. 

\begin{theorem}
\label{thm:main-queries}
 For all $\delta\in (0,1]$ and all $q \in [0,1]$, there exists $n_0\in \N,\eps_0 >0$ such that for every $\eps \in (0, \eps_0)$ and any $n \geq n_0$ the following holds. Any adaptive algorithm which reveals a path of length $\ell \geq \frac{\delta}{\eps}$ with probability at least $q$ in $G \sim G(n,\frac{1+\eps}{n})$ requires $\Omega(\frac{q\ell \delta}{p \eps})$ queries. 
\end{theorem}

We prove \Cref{thm:main-queries} by making the following improvement on a technical result that is the main ingredient in the proof of \Cref{thm:FKSV} in \cite{ferber-2017}.

\begin{theorem}
    \label{thm:CoverDisjointPaths}
    There exists $\eps_0>0$ such that for all $\delta \in (0,1]$ and $\eps \in (0, \eps_0)$, with high probability, the following holds.
    Let $S$ be a set of vertex-disjoint paths in $G \sim G(n,\frac{1+\eps}{n})$, each of length at least $\frac{\delta}{\eps}$.
    Then $S$ covers $O(\eps^2 \delta^{-2} n)$ vertices.
\end{theorem}

For a graph $G$ and an integer $\ell$, let $\cov_\ell(G)$ denote the maximum number of vertices that can be covered by a system of vertex-disjoint paths of length at least $\ell$ in $G$. In \cite{ferber-2017}, it was observed that to estimate $\cov_\ell(G)$ when $G\sim G(n,\frac{1+\eps}{n})$ it is actually sufficient to bound $\cov_\ell(T)$ for a random labelled tree $T$ of any given size $t$. Somewhat surprisingly, the insight that $\cov_\ell(T)=0$ when $T$ does not contain any path of length $\ell$, together with the trivial bound $\cov_\ell(T)\leq t$ for all other trees, is ultimately sufficient to prove \Cref{thm:FKSV}.

The key ingredient of our improvement is a more sophisticated bound on $\cov_\ell(T)$ for a typical labelled tree $T$. We obtain this bound through \Cref{lem:Counting-lemma} by counting the number of pairs ($e$,$T$) such that $e$ is an edge that lies near the centre of a long path of size at least $\ell$ of a random labelled tree on $t$ vertices.

The deduction of \Cref{thm:main-queries} from \Cref{thm:CoverDisjointPaths} proceeds in the same manner as that of \Cref{thm:FKSV} from the result in \cite{ferber-2017} that is analogous to \Cref{thm:CoverDisjointPaths}. For completeness, we include this proof in Appendix~\ref{ap:proof}.

Having identified the number of paths of length at least $\ell$ that can be packed into a random labelled tree as useful information for progress on \Cref{thm:FKSV}, the authors of \cite{ferber-2017} asked the following question, which we believe is of independent interest.

\begin{question}\label{question:tree-cover}
    Given $a=a(t) \in \mathbb{N}$ and $b=b(t) \in \mathbb{N}$, what is the probability that a random labelled tree on $t$ vertices contains $b$ vertex-disjoint paths, each of length at least $a$?
\end{question}

Because each path of length longer than $2\ell$ may be split into paths of length in $[\ell,2\ell]$, the maximum number of vertex-disjoint paths of length at least $\ell$ in a tree $T$ lies between $\cov_\ell(T)/2\ell$ and $\cov_\ell(T)/\ell$. Thus, by our aforementioned result that $\EE[\cov_\ell(T)]=\Theta(t/\ell)$, the average number of paths of length at least $\ell$ that can be packed into $T$ is $\Theta(t/\ell^2)$, as long as $\ell=O(\sqrt{t})$. What is more, if $\ell = o((t/\log t)^{1/3})$, we can use Talagrand's inequality to show that $\cov_\ell(T)$ is tightly concentrated around its mean, thus giving a partial answer to \Cref{question:tree-cover}. For the avoidance of doubt, the following theorem is phrased in terms of (unlabelled) trees on the fixed vertex set $[t]$ which is conceptually equivalent to labelled trees on $t$ unspecified vertices.

\begin{theorem}
\label{thm:ConcentrationPathPackingTree}
    There exist $C_1,C_2>0$ such that for all $t\in \N$ the following hold. Let $T$ be sampled uniformly at random from all trees on vertex set $[t]$. For all integers $\ell\in [1,t]$, we have $\EE[\cov_\ell(T)]\leq \frac{C_1t}{\ell}$, and for all integers $\ell\in [1,C_2t]$, we have $\EE[\cov_\ell(T)]\geq C_2t/\ell$.
    
    Furthermore, there exists $c>0$ such that for all integers $\ell \leq c\sqrt{t}$ and $\delta\in (0,1)$, if $T$ is sampled uniformly at random from all trees on vertex set $[t]$, then we have
    \begin{align*}
        \prob_T\left(\vert \cov_\ell(T) - \EE[\cov_\ell(T)]\vert > \frac{\delta t}{\ell}\right)\leq \frac{ct^2}{\delta e^{c\delta^2t/\ell^3}}.
    \end{align*}
\end{theorem}

The remainder of the paper is structured as follows. In \Cref{sec:prelim} we introduce some notation and tools, most of them already appearing in \cite{ferber-2017}. In \Cref{sec:MainCountingLemma} we state and prove the main counting lemma which is the crux of our proofs. In \Cref{sec:proofs} we proceed to the proof of \Cref{thm:CoverDisjointPaths}. In \Cref{sec:paths-in-tree} we prove \Cref{thm:ConcentrationPathPackingTree}.

\section{Preliminaries}
\label{sec:prelim}
For a positive integer $n$, let $[n] = \{1, \ldots, n\}$, and for a set $S$, let $S^{(n)}$ denote the set of all subsets of $S$ of size $n$. For functions $f$ and $g$, we denote $f \ll g$ if and only if there exist constants $N, k > 0$ such that for every $x > N$, we have $f(x) < k g(x)$. Note that this is different from the usual use of $\ll$ when studying random graphs.

We record the following version of Stirling's formula for later use.
\begin{lemma}
\label{lem:Stirling}
    There exists two absolute positive constants $K_1$ and $K_2$ such that for every integer $N$, we have $K_1 \sqrt{N} (\frac{N}{e})^N \geq N! \geq K_2 \sqrt{N} (\frac{N}{e})^N$.
\end{lemma}

We will make use of the following concentration inequality for the edge exposure martingale in $G(n,p)$ which is an easy consequence of \cite[Theorem 7.4.3]{alon2016probabilistic}.

\begin{lemma}
\label{lem:Martingale-exposure}
    Let $X$ be a random variable in the probability space $G(n,p)$ such that we have $|X(H_1)-X(H_2)| \leq C$ if $H_1$ and $H_2$ differ in at most one edge. Then, for any positive $\alpha < 2 \sqrt{n^2p}$, we have 
    \begin{align*}
        \prob(|X-\EE[X]| > C \alpha \sqrt{n^2p}) \leq 2e^{-\frac{\alpha^2}{4}}
    \end{align*}
\end{lemma}

The \emph{2-core} of a graph $G$ is the maximal induced subgraph of G in which all vertices have degree at least 2. The next lemma is a simple consequence of Theorem 5.4 of \cite{janson2011random} and Theorem 3 of \cite{pittel1990tree}.

\begin{lemma}
\label{lem:RandomGraphsProperties}
    Let $\eps>0$ be fixed, and let $p=\frac{1+\eps}{n}$. Then there exists $\eps_0>0$ such that for every $\eps \in (0, \eps_0)$, the following three properties hold whp in $G \sim G(n,p)$.
    \begin{enumerate}[(a)]
        \item the largest connected component of $G$ has size between $\eps n$ and $3 \eps n$,
        \item the second largest connected component of $G$ has size at most $20 \log n/\eps^2$,
        \item the $2$-core of the largest connected component of $G$ has size at most $2\eps^2 n$.
    \end{enumerate}
\end{lemma}

For $n \geq 1$, a path $P_n$ is a graph on $n$ vertices $v_1, \ldots, v_n$ with edges $v_1 v_2, \ldots, v_{n-1}v_n$. 

A \emph{rooted tree} is a connected graph with no cycles, where one of the vertices of $T$ has been marked as the \emph{root}. Recall that Cayley's formula\footnote{Cayley's formula is often phrased as counting the number of \emph{labelled trees} on $t$ vertices.} states that the number of trees on vertex set $[t]$ is $t^{t-2}$. It follows that the number of rooted trees on vertex set $[t]$ vertices is $t^{t-1}$. 

If $T$ is a tree rooted at $v$, then the \emph{depth} $\depth(w)$ of a vertex $w$ is the distance from $w$ to $v$ in $T$. The \emph{height} of $T$, denoted by $h(T)$ is the maximum depth of a vertex in $T$, while the \emph{width} of $T$ is given by
\begin{equation*}
    w(T)=\max_{k\in [h(T)]}\vert \{w\in V(T):\depth(w)=k\} \vert.
\end{equation*}
If $u$ is a vertex in a tree $T$ without a root, we refer to the height of $T$ rooted at $u$ as the \emph{height of $u$}.

A \emph{Galton-Watson tree} is a rooted tree which is constructed recursively at random. Starting from the root, each vertex is given a random number of children independently of the other vertices according to a distribution $\xi$ on $\mathbb{N}_0$. The distribution $\xi$ is called the \emph{offspring distribution} and the tree is referred to as the $\xi$-Galton-Watson tree. For concreteness, we say that the vertex set of a Galton-Watson tree $T$ is the integer interval $[\vert T\vert]$.

In \cite{addario2013sub}, Addario-Berry, Devroye, and Janson established sub-gaussian tail bounds for the distributions of the width and height of a Galton-Watson tree conditioned on the number of its vertices. In the following two theorems, $\xi$ is a distribution on $\N_0$ with expectation 1 and positive, but finite, variance.

\begin{theorem}\label{thm:height-bound}
There exist $c_1, C_1>0$ such that for a $\xi$-Galton-Watson tree $T$ and all $x,t\in \N$, 
\begin{equation*}
    \prob(h(T) \geq x \mid \vert V(T)\vert=t) \leq C_1e^{-c_1 x^2/t}.
\end{equation*}
\end{theorem}

\begin{theorem}\label{thm:width-bound}
There exist $c_2, C_2>0$ such that for a $\xi$-Galton-Watson tree $T$ and all $x,t\in \N$,
\begin{equation*}
    \prob(w(T) \geq x\mid \vert V(T)\vert=t) \leq C_2e^{-c_2 x^2/t}.
\end{equation*}
\end{theorem}

The relevance of these theorems for our paper stems from the following theorem from \cite[Section 6.4]{devroye1998branching}.

\begin{theorem}\label{thm:gw-uniform}
    If $\xi$ is the Poisson(1) distribution and $T$ is a $\xi$-Galton-Watson tree, then the distribution of $T$ conditioned on the event $\vert V(T)\vert = t$ is the same as the uniform distribution over rooted trees on vertex set $[t]$.
\end{theorem}

The following lemma is a corollary of the previous three results, and it will be central in our proof of \Cref{thm:CoverDisjointPaths}.

\begin{lemma}
\label{lem:HeightLabelledTree}
    There exist constants $c_1, c_2, C_1,C_2 >0$ such that for all $x,t\in \N$, if $T$ is sampled uniformly at random from all rooted trees on $[t]$, then 
    \begin{equation}\label{eq:height-upper-bound}
        \prob(h(T) \geq x) \leq C_1e^{-c_1 x^2/t},
    \end{equation}
    and 
    \begin{equation}\label{eq:height-lower-bound}
        \prob(h(T) \leq t/x) \leq C_2e^{-c_2 x^2/t}.
    \end{equation}
\end{lemma}

\begin{proof}
    The first equation \eqref{eq:height-upper-bound} follows immediately from Theorems \ref{thm:height-bound} and \ref{thm:gw-uniform}, while \eqref{eq:height-lower-bound} follows from Theorems \ref{thm:width-bound}, \ref{thm:gw-uniform} and the additional observation that $h(T)w(T)\geq t$.
\end{proof}

We need the following well-known fact about Poisson($\mu$)-Galton-Watson tree (see for instance Section 6 of \cite{devroye1998branching}).

\begin{lemma}
\label{lem:distrib-Poisson-GW}
    Let $\mathcal{T}$ be a Poisson($\mu$)-Galton-Watson tree. Then
    $$\prob(|\mathcal{T}|=t)=\frac{t^{t-1}(\mu e^{-\mu})^t}{\mu t!}.$$
\end{lemma}

Ding, Lubetzky and Peres \cite{ding-2014} gave a full characterisation of the structure of the giant component of $G \sim G(n,p)$ in the strictly supercritical regime, i.e. where $p=\frac{1+\eps}{n}$ for some constant $\eps>0$. We will use the following consequence of their work, which was already used in a slightly different form in \cite{ferber-2017}.

\begin{lemma}
\label{lem:consequence-Ding-anatomy} 
    Let $p=\frac{1+\eps}{n}$ for some constant $\eps>0$.
    Let $\mathcal{C}_1$ be the largest connected component of $G \sim G(n,p)$, let $\mathcal{C}^{(2)}_1$ be the $2$-core of $\mathcal{C}_1$, and let $\mathcal{C}_1 \setminus \mathcal{C}^{(2)}_1$ be the graph obtained from $\mathcal{C}_1$ by deleting all edges from $\mathcal{C}^{(2)}_1$.
    Let $0 < \mu <1$ be such that $\mu e^{-\mu}=(1+\eps)e^{-(1+\eps)}$ and consider $2 \eps^2 n$ independent Poisson($\mu$)-Galton-Watson trees $T_1, \dots, T_{2\eps^2 n}$.
    Then, for every $\ell$ and $M$, if whp the disjoint union of $T_1, \dots, T_{2 \eps^2 n}$ does not contain a set of vertex-disjoint paths of length $\ell$ covering at least $M$ edges, then the same holds whp for $\mathcal{C}_1 \setminus \mathcal{C}^{(2)}_1$.
\end{lemma}

Finally, to prove the concentration part of \Cref{thm:ConcentrationPathPackingTree}, we will also require Talagrand's inequality. Let $\Omega=\prod_{i=1}^n\Omega_i$ be a probability space where each $\Omega_i$ is a finite probability space and $\Omega$ has the product measure. A random variable $X\colon \Omega\rightarrow \R$ is called \emph{$K$-Lipschitz} if $\vert X(\omega)-X(\omega')\vert \leq K$ whenever $\omega,\omega'\in \Omega$ differ in at most one coordinate. Let $f\colon \N\rightarrow \N$. We say that $X$ is \emph{$f$-certifiable} if for all $\omega\in \Omega$ and $s\in \N$, if $X(\omega)\geq s$ there exists
$I \subset [n]$ with $\vert I\vert\leq f(s)$ so that all $\omega' \in \Omega$ that agree with $\omega$ on $I$ satisfy $X(\omega')\geq s$. The following theorem is a simple corollary to Talagrand's inequality (see Theorem 7.7.1 of \cite{alon2016probabilistic}).

\begin{theorem}\label{thm:talagrand}
    Let $f\colon \N\rightarrow \N$ be a function and let $X\colon \Omega\rightarrow \R$ be $K$-Lipschitz and $f$-certifiable. Then for all $\mu, \alpha\in \R$,
    \begin{align*}
        \prob\left(X\leq \mu - \alpha K\sqrt{f(\mu)}\right)\prob(X\geq \mu) \leq e^{-\alpha^2/4}.
    \end{align*}
\end{theorem}

\section{The main counting lemma}
\label{sec:MainCountingLemma}

For a graph $G$, we say that an edge $uv$ of $G$ is $m$-centred if there exist two vertex-disjoint paths of length at least $m-1$ in $G - uv$ starting at $u$ and $v$, respectively. For $t\geq m$, we define $\midpath(t,m)$ as the set of pairs $(uv,T)$ such that $T$ is a tree on $[t]$ and $uv$ is an $m$-centred edge of $T$.

Recall that for a graph $G$ and an integer $\ell$, $\cov_\ell(G)$ is the maximum number of vertices that can be covered by a system of vertex-disjoint paths of length at least $\ell$ in $G$. The key insight of our proofs is that the number of $m$-centred edges in a tree $T$ can be used to estimate $\cov_\ell(T)$. Indeed, if $v_0\ldots v_{k}$ is a path in $T$ for $k\geq 3m$ and $m-1 \leq i \leq k-m$, then $T-v_iv_{i+1}$ contains two vertex-disjoint paths of length $m-1$ that start at $v_i$ and $v_{i+1}$, respectively. In other words, $v_iv_{i+1}$ is $m$-centred in $T$. Since the interval $[m-1,k-m]$ has length $k-2m+1>\frac{k+1}{3}$, it follows that
\begin{equation}\label{eq:cov-leq-centred}
    \cov_{3m}(T) \leq 3\vert \{\text{$m$-centred edges in $T$}\}\vert.
\end{equation}

Conversely, we may choose an arbitrary root $r$ in $T$ and construct a system of vertex-disjoint paths in $T$ greedily by taking the longest monotone (in the tree order) path that is disjoint from all previously included paths until no path of length $\ell$ is left. The resulting systems of paths will include all vertices that are incident to $\ell$-centred edges in $T$. Thus, we have
\begin{equation}\label{eq:centred-leq-cov}
    \vert \{\text{$\ell$-centred edges in $T$}\}\vert \leq \cov_\ell(T).
\end{equation}

The following lemma will allow us to estimate the number of $m$-centred edges in a typical tree.

\begin{lemma}
\label{lem:Counting-lemma}
    There exist constants $c_1, C_1>0$ such that for all $t\in \N$ and $m<t$,
    \begin{align*}
        |\midpath(t,m)| < C_1\sum_{k=m}^{t/2} t^{t-1}k^{-3/2} e^{-c_1\frac{m^2}{k}}.
    \end{align*}
    Furthermore, there exist constants $c_2, C_2>0$ such that
    \begin{align*}
       |\midpath(t,m)| > C_2\sum_{k=c_2 m^2}^{t/2} t^{t-1}k^{-3/2}.
    \end{align*}
\end{lemma}

\begin{proof}
    For $m,t\in \N$, $S\subset [t]$ and $u\in S$ let $\deep_m(S;u)$ be the set of trees on $S$ rooted at $u$ whose height is at least $m-1$. For $v\in [t]\setminus S$ and a rooted tree $T_v\in \deep_m(S,v)$, we write $T_u+T_v$ for the tree on $[t]$ with edge set $E(T_u)\cup E(T_v)\cup \{uv\}$. Note that $uv$ is an $m$-centred edge of $T_u+T_v$ so that $(uv,T_u+T_v)\in \midpath(t,m)$. Thus the following map is well-defined.
    \begin{align*}
        \phi\colon  \bigcup_{\substack{S\subset [t]\\ u\in S\\ v\in [t]\setminus S}} \deep_m(S;u)\times \deep_m([t]\setminus S;v) &\rightarrow \midpath(t,m)  & (T_u,T_v) &\mapsto (uv,T_u + T_v).
    \end{align*}
    For each $(uv,T)\in \midpath(t,m)$, if $T_u$ and $T_v$ are the components of $T-uv$ rooted at $u$ and $v$, respectively, then $\phi^{-1}(uv,T)=\{(T_u,T_v),(T_v,T_u)\}$. Thus, the domain of $\phi$ is exactly twice as large as its codomain $\midpath(t,m)$, and since $\deep_m(S;u)$ is empty unless $\vert S\vert \geq m$, we have
    \begin{align}\label{eq:midpath-bound1}
        \vert \midpath(t,m)\vert &= \frac{1}{2}\sum_{k=m}^{t-m} \sum_{\substack{S\in [t]^{(k)}\\ u\in S\\ v\in [t]\setminus S}} \vert \deep_m(S;u)\vert \cdot  \vert \deep_m([t]\setminus S;v)\vert\nonumber\\
        &=\frac{1}{2}\sum_{k=m}^{t-m} \binom{t}{k}k(t-k) \vert \deep_m([k];1)\vert \cdot  \vert \deep_m([t-k];1)\vert.
    \end{align}
    To bound this quantity from above, note that since there are exactly $k^{k-2}$ trees on $[k]$ overall, we may infer from \Cref{lem:HeightLabelledTree} that
    \begin{align*}
        \vert \deep_m([k];1)\vert \ll k^{k-2}e^{-cm^2/k}.
    \end{align*}
    Inserting this into \eqref{eq:midpath-bound1}, we find that
    \begin{align*}
        \vert \midpath(t,m)\vert &\ll \sum_{k=m}^{t-m} \binom{t}{k}k^{k-1}(t-k)^{t-k-1} e^{-cm^2(1/k+1/(t-k))}\\
        &\ll \sum_{k=m}^{t/2} \binom{t}{k}k^{k-1}(t-k)^{t-k-1} e^{-cm^2/k}\\
        &=\sum_{k=m}^{t/2} \frac{t!k^{k-1}(t-k)^{t-k-1}}{k!(t-k)!} e^{-cm^2/k}.
    \end{align*}
    By Stirling's formula (\Cref{lem:Stirling}), this can be bounded as
    \begin{align*}
        \vert \midpath(t,m)\vert &\ll \sum_{k=m}^{t/2} \frac{t^{t+1/2}}{k^{3/2}(t-k)^{3/2}} e^{-cm^2/k}\ll t^{t-1}\sum_{k=m}^{t/2} k^{-3/2} e^{-cm^2/k},
    \end{align*}
    as desired.

    To bound $\vert \midpath(t,m)\vert$ from below, we observe that by \Cref{lem:HeightLabelledTree},
    \begin{align*}
        \vert \deep_m([k];1) \vert \geq k^{k-2} (1-Ce^{-ck/m^2}).
    \end{align*}
    In particular, there exists a constant $c_2$ such that $\vert \deep_m([k];1) \vert \geq k^{k-2}/2$ when $k>c_2 m^2$. Therefore, we deduce from \eqref{eq:midpath-bound1} that
    \begin{align*}
        \vert \midpath(t,m)\vert \gg \sum_{k=c_2 m^2}^{t/2} \binom{t}{k}k^{k-1}(t-k)^{t-k-1}.
    \end{align*}
    As before, it follows from Stirling's formula that
    \begin{align*}
        \vert \midpath(t,m)\vert &\gg  t^{t-1}\sum_{k=c_2 m^2}^{t/2} k^{-3/2}. \qedhere
    \end{align*}
\end{proof}

\section{Proof of \Cref{thm:CoverDisjointPaths}}
\label{sec:proofs}

\begin{proof}[Proof of \Cref{thm:CoverDisjointPaths}]
    For every $\eps>0$, there exists a unique $\mu=\mu(\eps)\in (0,1)$ such that $\mu e^{-\mu}=(1+\eps)e^{-(1+\eps)}$. We choose $\eps_0$ such that for all $\eps<\eps_0$, we have $1-\mu(\eps)>\eps/2$ and $1+\eps \leq e^{\eps - \frac{\eps^2}{3}}$. Let $\delta\in (0,1]$ and $\eps\in (0,\eps_0)$. 
    
    As a path of length $m$ consists of $m-1$ edges and $m$ vertices, it is equivalent to show that whp any set of vertex-disjoint paths of lengths at least $\ell=\frac{\delta}{\eps}$ covers at most $D\eps^2 \delta^{-2} n$ edges, for some absolute constant $D > 0$. Note that the claim is trivial if $\ell<100$ so that we may assume that $\ell\geq 100$. Here and throughout the proof, we ignore rounding unless it is critical to the correctness of the proof.
    
    Let $\mathcal{C}_1$ be the largest connected component of $G$, let $\mathcal{C}^{(2)}_1$ be the $2$-core of $\mathcal{C}_1$, and let $\mathcal{C}_1 \setminus \mathcal{C}^{(2)}_1$ be the graph obtained from $\mathcal{C}_1$ by deleting the edges of $\mathcal{C}^{(2)}_1$.
    Furthermore, we define the following random variables:
    \begin{itemize}
        \item Let $X_{\ell}$ be the maximum number of edges covered by vertex-disjoint paths of length at least $\ell$ in components of $G$ of size at most $20 \log n/\eps^2$.
        \item Let $Y_{\ell}$ be the maximum number of edges covered by vertex-disjoint paths of length at least $\ell$ in $\mathcal{C}_1$.
        \item Let $Z_{\ell}$ be the maximum number of edges covered by vertex-disjoint paths of length at least $\ell/3$ in $\mathcal{C}_1 \setminus \mathcal{C}^{(2)}_1$.
    \end{itemize}
    By \Cref{lem:RandomGraphsProperties}(b), we have that whp the number of edges of $G$ covered by vertex-disjoint paths of length at least $\ell$ is at most $X_{\ell}+Y_{\ell}$.
    As in \cite{ferber-2017}, a simple counting argument shows that $Y_{\ell} \leq 6|\mathcal{C}^{(2)}_1|+6Z_{\ell}$, and hence by \Cref{lem:RandomGraphsProperties}(c), we have whp $Y_{\ell} \leq 12 \eps^2 n+6Z_{\ell}$. Therefore, to finish the proof, it suffices to prove that whp we have $X_{\ell} \ll \eps^2\delta^{-2} n$ and $Z_\ell \ll \eps^2\delta^{-2} n$.
    We start by proving that $X_{\ell} \ll \eps^2 \delta^{-2} n$ whp.

    Let $S\subset [n]$ and set $m=\ell/3$. 
    We define a random variable $X_S$ whose value depends on a case distinction. 
    If $G[S]$ is connected and there are no edges between $S$ and $[n]\setminus S$, then $X_S$ is defined as the number of edges in $G[S]$ that are $m$-centred. If either condition is not satisfied, we define $X_S$ to be 0. Observe next that if $G[S]$ is connected and $e$ is an $m$-centred edge in $G[S]$, then $G[S]$ has a spanning tree $T$ such that $e$ is also $m$-centred in $T$.
    Therefore, we have
    \begin{align}
    \label{eq:BoundFixedSTriples}
        \mathbb{E}[X_S] \leq |\mathcal{M}(t,m)|p^{t-1} (1-p)^{t(n-t)}.
    \end{align}
    Indeed, $|\mathcal{M}(t,m)|$ counts the number of pairs $(e,T)$ such that $T$ is a spanning tree on $S$ and $e$ is $m$-centred in $T$, $p^{t-1}$ is the probability that $T\subset G[S]$, and $(1-p)^{t(n-t)}$ is the probability that there are no edges between $S$ and $[n]\setminus S$.
    Since $X_\ell$ is bounded by the sum of all $X_S$ for which $S\subset [n]$ has size between $m$ and $20\log n/\eps^2$, inserting the upper bound from \Cref{lem:Counting-lemma} into \eqref{eq:BoundFixedSTriples}, yields
    \begin{align*}
         \EE[X_{\ell}] &\ll \sum_{t=m}^{\frac{20}{\eps^2} \log n} \binom{n}{t} \sum_{k=m}^{t/2} t^{t-1}k^{-3/2} e^{-c\frac{m^2}{k}} p^{t-1} (1-p)^{t(n-t)} \\
         &\ll \sum_{t=m}^{\frac{20}{\eps^2} \log n} \frac{n!}{t!(n-t)!} \sum_{k=m}^{t/2} t^{t-1}k^{-3/2} e^{-c\frac{m^2}{k}} \frac{e^{\eps(t-1)}}{n^{t-1}} e^{-\frac{1+\eps}{n}t(n-t)}.
    \end{align*}
    Simplifying and applying \Cref{lem:Stirling}, we obtain
    \begin{align}
    \label{eq:proof1}
        \EE[X_{\ell}] &\ll n \sum_{t=m}^{\frac{20}{\eps^2} \log n} \frac{n^{n-t+1/2}}{t^{t+1/2}(n-t)^{n-t+1/2}} \sum_{k=m}^{t/2} t^{t-1}k^{-3/2} e^{-c\frac{m^2}{k}} e^{-t}.
    \end{align}
    Moreover, we have 
    \begin{align}
    \label{eq:proof2}
        \frac{n^{n-t+1/2}}{(n-t)^{n-t+1/2}} = \frac{1}{(1-t/n)^{n-t+1/2}} = e^{t+O(\frac{t^2}{n})} \ll e^{t}.
    \end{align}
    Inserting \eqref{eq:proof2} into \eqref{eq:proof1}, we obtain
    \begin{align*}
        \EE[X_{\ell}] &\ll n \sum_{t=m}^{\frac{20}{\eps^2} \log n} \sum_{k=m}^{t/2} k^{-3/2}t^{-3/2} e^{-c\frac{m^2}{k}} \\
        &\ll n \sum_{t=m}^{\infty} \sum_{k=m}^{t/2} k^{-3/2}t^{-3/2} e^{-c\frac{m^2}{k}} \\
        &\ll n \sum_{t=m}^{\infty} t^{-3/2} \int_{m}^{t/2} k^{-3/2}e^{-c\frac{m^2}{k}} \diff k,
    \end{align*}
    where the integral comparison that yields the last inequality is valid because the integrand is positive and for all $x,y\geq m \geq 1$ such that $\vert x-y\vert \leq 1$, we have
    \begin{equation*}
        \frac{x^{-3/2}e^{-c\frac{m^2}{x}}}{y^{-3/2}e^{-c\frac{m^2}{y}}}=\left(\frac{y}{x}\right)^{3/2}e^{c\left( \frac{m^2}{y}-\frac{m^2}{x} \right)}\leq 2^{3/2} e^{c}.
    \end{equation*}
Performing the substitution $k={m}^2/x^2$, we obtain
    \begin{align*}
     \EE[X_{\ell}] &\ll n{m}^{-1} \sum_{t=m}^{\infty} t^{-3/2} \int_{\sqrt{\frac{2}{t}}m}^{\sqrt{m}} e^{-cx^2} \diff x \\
     &\ll n{m}^{-1} \sum_{t=m}^{\infty} t^{-3/2} \int_{0}^{\sqrt{m}} \indic{x \geq \sqrt{\frac{2}{t}}m} e^{-cx^2} \diff x\\
     &\ll n{m}^{-1} \int_{0}^{\sqrt{m}} e^{-cx^2}\sum_{t=\max(2m^2/x^2, m)}^{\infty} t^{-3/2}  \diff x\\
    &\ll n{m}^{-1} \int_{x=0}^{\sqrt{m}} e^{-cx^2} \int_{t=2m^2/x^2}^{\infty} t^{-3/2} \diff t \diff x \\
     &\ll n{m}^{-1} \int_{x=0}^{\sqrt{m}} e^{-cx^2} (2m^2/x^2)^{-1/2} \diff x \\
     &\ll n{m}^{-2} \int_{0}^{\infty} xe^{-cx^2} \diff x.
    \end{align*}
    Since the moments of a Gaussian random variable exist and are finite, we may insert the value of $m$ to obtain $\EE[X_{\ell}] \ll n \eps^2\delta^{-2}$. An application of \Cref{lem:Martingale-exposure} gives the desired conclusion. \\

    We now move on to bounding $Z_\ell$. Let $0 < \mu <1$ be such that $\mu e^{-\mu}=(1+\eps)e^{-(1+\eps)}$ and consider $2 \eps^2 n$ independent Poisson($\mu$)-Galton-Watson trees $T_1, \dots, T_{2\eps^2 n}$. By \Cref{lem:consequence-Ding-anatomy}, to prove that $Z_\ell \ll \eps^2 \delta^{-2} n$ whp, it is enough to show that $M_\ell$, the number of edges that can be covered by vertex-disjoint paths in $T_1\cup \ldots \cup T_{2\eps^2 n}$ of length at least $\ell/3$, satisfies $M_\ell\ll \eps^2 \delta^{-2} n$ whp.

    To bound $M_\ell$, we employ the second moment method, as was done in \cite{ferber-2017}. However, \Cref{lem:Counting-lemma} will allow us to obtain better bounds once more.
    Set $m=\ell/9$ and let $T$ be a Poisson($\mu$)-Galton-Watson tree and let $S_\ell$ be the number of $m$-centred edges in $T$. Thus, $\EE[M_\ell]\leq 6\eps^2n \EE[S_\ell]$. It is clear that
    \begin{align*}
        \EE[S_{\ell}] &\leq 3\sum_{t \geq m} \mathbb{P}(|T|=t)|\frac{|\mathcal{M}(t,m)|}{t^{t-2}}.
    \end{align*}
    Inserting \Cref{lem:distrib-Poisson-GW,lem:Counting-lemma} into the above, we obtain
    \begin{align*}
        \EE[S_{\ell}] &\ll \sum_{t \geq m} \frac{1}{t^{t-2}} \frac{t^{t-1}(\mu e^{-\mu})^t}{\mu t!} \sum_{k=m}^{t/2} t^{t-1}k^{-3/2} e^{-c\frac{m^2}{k}} \\
        &\ll \sum_{t \geq m} \frac{t (1+\eps)^te^{-(1+\eps)t}}{t!} \sum_{k=m}^{t/2} t^{t-1}k^{-3/2} e^{-c\frac{m^2}{k}}.
    \end{align*}
    Using that $(1+\eps)^t \leq e^{\eps t- \frac{\eps^2}{3}t}$ and \Cref{lem:Stirling}, we have
    \begin{align*}
        \EE[S_{\ell}] &\ll \sum_{t \geq m} \sum_{k=m}^{t/2} k^{-3/2}t^{-1/2} e^{-\frac{\eps^2}{3}t} e^{-c\frac{m^2}{k}} \\
        &\ll \sum_{t \geq m} e^{-\frac{\eps^2}{3}t} t^{-1/2} \int_{k=m}^{t/2}  k^{-3/2}e^{-c\frac{m^2}{k}}\diff k.
    \end{align*}
    As for the computation of $\EE[X_{\ell}]$, we substitute $k={m}^2/x^2$ in the integral to obtain
    \begin{align*}
        \EE[S_{\ell}] &\ll {m}^{-1} \sum_{t=m}^{\infty} e^{-\frac{\eps^2}{3}t} t^{-1/2} \int_{\sqrt{\frac{2}{t}}m}^{\sqrt{m}} e^{-cx^2} \diff x \\
        &\ll {m}^{-1} \sum_{t=m}^{\infty} \int_{x=\sqrt{\frac{2}{t}}m}^{\sqrt{m}} \indic{x \geq \sqrt{\frac{2}{t}}m} e^{-\frac{\eps^2}{3}t} t^{-1/2} e^{-cx^2} \diff x \\
        &\ll {m}^{-1} \int_{x=0}^{\sqrt{m}} e^{-cx^2} \sum_{t=2{m}^2/x^2}^{\infty}  e^{-\frac{\eps^2}{3}t} t^{-1/2} \diff x \\
        &\ll {m}^{-1} \int_{x=0}^{\sqrt{m}} e^{-cx^2} \int_{t=2{m}^2/x^2}^{\infty}  e^{-\frac{\eps^2}{3}t} t^{-1/2} \diff t \diff x.
    \end{align*}
    Performing the second substitution $t=y^2/\eps^2$ in the integral, we obtain
    \begin{align*}
        \EE[S_{\ell}] &\ll {m}^{-1} \int_{x=0}^{\sqrt{m}} e^{-cx^2} \int_{y=\sqrt{2}m \eps/x}^{\infty}  \frac{e^{-\frac{y^2}{3}}}{\eps}  \diff y \diff x \\
        &\ll \delta^{-1} \int_{x=0}^{\infty} e^{-cx^2} \int_{y=0}^{\infty}  e^{-\frac{y^2}{3}} \diff y \diff x \\
        &\ll \delta^{-1},
    \end{align*}
    showing that $\EE[M_\ell]\ll \eps^2\delta^{-1}n$.
    
    For the second moment, we use the same computation as \cite{ferber-2017}. Using \Cref{lem:distrib-Poisson-GW}, which specifies the distribution of $\vert T\vert$, we can calculate that $\EE[\vert T\vert^2]=(1-\mu)^{-3}$. Thus,
    \begin{align*}
        \Var[M_{\ell}] \leq 18\eps^2 n \EE[S^2_{\ell}] \leq 18\eps^2 n\EE[|T|^2] =\frac{18\eps^2 n}{(1-\mu)^3} \leq \frac{144n}{\eps},
    \end{align*}
    where the last inequality follows from $1-\mu \geq \frac{\eps}{2}$.
    Since $\EE[M_\ell]\ll \eps^2\delta^{-1}n$ and $\Var[M_{\ell}] \ll n\eps^{-1}$ as well, it follows from Chebyshev's inequality that $M_{\ell} \ll \eps^2 \delta^{-1} n$ whp, as desired.
\end{proof}

\section{Proof of \Cref{thm:ConcentrationPathPackingTree}}
\label{sec:paths-in-tree}

Recall that $\cov_\ell(T)$ is the maximum number of vertices in $T$ that can be covered by vertex-disjoint paths of length at least $\ell$. We begin by estimating the expected value of $\cov_\ell(T)$.

\begin{proposition}\label{prop:path-packing-mean}
There exists a constant $C$ such that for all integers $1 \leq \ell\leq t$, we have $\EE[\cov_\ell(T)]\leq \frac{Ct}{\ell}$, where $T$ is sampled uniformly at random from all trees on vertex set $[t]$. Furthermore, there exists a constant $c>0$ such that for all $t\in \N$ and $1\leq \ell<c\sqrt{t}$, we have $\EE[\cov_\ell(T)]\geq ct/\ell$.
\end{proposition}

\begin{proof}
    For a tree $T$ on $[t]$ and $m\leq t$, let $X_m(T)$ be the number of pairs $(e,T)$ in $\midpath(t,m)$, i.e., the number of $m$-centred edges in $T$. Recall that by \eqref{eq:cov-leq-centred}, we have $\cov_{3m}(T)\leq 3X_m(T)$, while by \eqref{eq:centred-leq-cov}, we have $X_\ell(T)\leq \cov_\ell(T)$. Observe further that since $\cov_{\ell'}(T)\leq \cov_\ell(T)$ for $\ell'>\ell$, it is enough to show the upper bound in the claim when $\ell$ is divisible by three. Thus, to prove the claim, it suffices to show that $\EE[X_m]=O(t/m)$ for $1\leq m\leq t$, and that there exists $c>0$ such that $\EE[X_m]=\Omega(t/m)$ for $1\leq m \leq c\sqrt{t}$.
    
    To bound $\EE[X_m]$ from above, we make similar calculations to those in the proof of \Cref{thm:CoverDisjointPaths}. By \Cref{lem:Counting-lemma}, we have
    \begin{align*}
        \EE [X_m(T)] &\ll t\sum_{k=m}^{t/2} k^{-3/2} e^{-c\frac{t^2}{k}} \ll t\int_m^{t/2} k^{-3/2} e^{-c\frac{t^2}{k}} \diff k \ll t\int_m^{t/2} k^{-3/2} e^{-c\frac{m^2}{k}} \diff k.
    \end{align*}
    By substituting $k=m^2/x^2$, we obtain
    \begin{align*}
        \EE [X_m(T)] &\ll \frac{t}{m}\int_{m/\sqrt{t}}^{\sqrt{m}} e^{-cx^2} \diff x\ll\frac{t}{m},
    \end{align*}
    as desired.

    It remains to bound $\EE[X_m(T)]$ from below. By the lower bound in \Cref{lem:Counting-lemma}, there exists a constant $c_2$ such that
    \begin{align*}
        \EE[X_m(T)] \gg t\sum_{k=c_2\ell^2}^{t/2} k^{-3/2}.
    \end{align*}
    Assuming that $c_2m^2<t/4$, this gives
    \begin{align*}
        \EE[X_m(T)] \gg t \int_{c_2 m^2}^{t/2} k^{-3/2} \diff k \gg \frac{t}{\sqrt{c_2 m^2}}-\frac{t}{\sqrt{t/2}} \gg \frac{t}{m},
    \end{align*}
    as desired.
\end{proof}

As a function of the edges of $T$, $\cov_\ell(T)$ is continuous in the sense that changing a single edge can change $\cov_\ell$ by at most $2\ell$. If $\ell$ is sufficiently small compared to $t$, we can use Talagrand's inequality (see \Cref{thm:talagrand}) to exploit this continuity to show that $\cov_\ell$ is concentrated around its mean.

\begin{proof}[Proof of \Cref{thm:ConcentrationPathPackingTree}]
    For $v\in [t-1]$ let $\Omega_v=[t]$ and consider the product space $\Omega=\prod_{v=1}^{t-1} \Omega_v$ with uniform probability measure $\prob_\omega$. Each element $\omega\in \Omega$ can be associated with a graph with vertex set $[t]$ and edge set $\{uv:\omega_u=v\lor \omega_v=u\}$. In particular, each tree $T$ on $[t]$ is associated with the unique element $\omega\in \Omega$ for which $\omega_v$ is the first vertex on the path from $v$ to $t$. For $\omega\in \Omega$, we define $\cov_\ell(\omega)$ as $\cov_\ell(G_\omega)$ for the unique graph $G_\omega$ associated with $\omega$.

    Observe that the random variable $\omega \mapsto \cov_\ell(G_\omega)$ is $2\ell$-Lipschitz. Furthermore, if $G_\omega$ has a system of vertex-disjoint paths of length at least $\ell$ that cover $m$ vertices, then the $m-1$ edges of this path system certify that $\cov_\ell(G_\omega)\geq m$. Hence, $\cov_\ell$ is $f$-certifiable for $f(s)=s-1$. Therefore, we may apply \Cref{thm:talagrand} to see that for any $\alpha,\mu >0$,
    \begin{align*}
        \prob_\omega\left(\cov_\ell(G_\omega)\leq \mu -   \alpha\ell\sqrt{\mu}\right) \prob_\omega(\cov_\ell(G_\omega)\geq \mu)\leq e^{-\alpha^2/4}.
    \end{align*}
    Since the size of $\Omega$ is $t^{t-1}$ and there are $t^{t-2}$ trees on $[t]$, the set $\{\omega:G_\omega \text{ is a tree}\}$ has density $1/t$ in $\Omega$, it follows that for the uniform probability measure $\prob_T$ over all trees on $[t]$, we have
    \begin{equation}\label{eq:packing-talagrand}
        \prob_T\left(\cov_\ell(T)\leq \mu -  \alpha\sqrt{\ell^2 \mu}\right) \prob_T(\cov_\ell(T)\geq \mu)\leq t^2e^{-\alpha^2/4}.
    \end{equation}
    Let $\delta\in (0,1)$. By \Cref{prop:path-packing-mean}, there exists $c>0$ such that $\EE_T[\cov_\ell(T)]\geq ct/\ell$, independently of $t$ as long as $\ell\ll c\sqrt{t}$. Thus, to bound $\prob_T(\cov_\ell(T)\geq (1+\delta)\EE_T[\cov_\ell(T)])$, we may apply \eqref{eq:packing-talagrand} with $\mu=(1+\delta)\EE_T[\cov_\ell(T)]$ and
    \begin{equation*}
        \alpha=\frac{\delta}{4}\sqrt{\frac{\EE_T[\cov_\ell(T)]}{\ell^2}}\geq \frac{\delta}{4}\sqrt{\frac{ct}{\ell^3}},
    \end{equation*}
    so that 
    \begin{equation*}
        \mu-\alpha\sqrt{\ell^2\mu}=(1+\delta)\EE_T[\cov_\ell(T)]-\frac{\delta}{4}\sqrt{(1+\delta)\EE_T[\cov_\ell(T)]^2}\geq (1+\delta/2)\EE_T[\cov_\ell(T)].
    \end{equation*}
    Inserting this into \eqref{eq:packing-talagrand}, we obtain
    \begin{equation*}
        \prob_T(\cov_\ell(T)\geq (1+\delta)\EE_T[\cov_\ell(T)])\leq \frac{t^2e^{-\frac{ct\delta^2}{64\ell^3}}}{\prob_T\left(\cov_\ell(T)\leq (1+\delta/2)\EE_T[\cov_\ell(T)]\right)}.
    \end{equation*}
    By Markov's inequality, the denominator is at least $1-1/(1+\delta/2)\geq \delta/3$ so that
    \begin{equation}\label{eq:concentration-upper}
        \prob_T(\cov_\ell(T)\geq (1+\delta)\EE_T[\cov_\ell(T)])\leq \frac{3t^2e^{-\frac{ct\delta^2}{64\ell^3}}}{\delta}.
    \end{equation}

    On the other hand, if we let $\mu=\EE_T[\cov_\ell(T)]$ and $\alpha=\delta\sqrt{ct/\ell^3}$ and denote $\prob_T(\cov_\ell(T)<(1-\delta)\EE_T[\cov_\ell(T)])$ by $q_\delta$, then
    \begin{equation*}
        \prob_T(\cov_\ell(T)\geq \EE_T[\cov_\ell(T)])\leq t^2q_\delta^{-1}e^{-\frac{ct\delta^2}{\ell^3}}.
    \end{equation*}
    Since $\cov_\ell(T)$ is bounded by $t$, $\EE_T[\cov_\ell(T)]$ can be bounded as
    \begin{align*}
        \EE_T[\cov_\ell(T)]&\leq q_\delta(1-\delta)\EE_T[\cov_\ell(T)]+(1-q_\delta)\EE_T[\cov_\ell(T)]+\prob_T(\cov_\ell(T) > \EE_T[\cov_\ell(T)])t\\
        &\leq (1-q_\delta\delta)\EE_T[\cov_\ell(T)] + t^3q_\delta^{-1}e^{-\frac{ct\delta^2}{\ell^3}}, 
    \end{align*}
    which implies that
    \begin{equation}\label{eq:concentration-lower}
        \prob_T(\cov_\ell(T)<(1-\delta)\EE_T[\cov_\ell(T)])=q_\delta \leq \left(\frac{t^3}{\delta\EE_T[\cov_\ell(T)]}e^{-\frac{ct\delta^2}{\ell^3}}\right)^{1/2}\leq t\sqrt{\frac{\ell}{c\delta}} e^{-\frac{ct\delta^2}{2\ell^3}}.
    \end{equation}
    We may now apply a union bound for \eqref{eq:concentration-upper} and \eqref{eq:concentration-lower}, and by renaming the constant $c$ appropriately, we obtain the claim.
\end{proof}

The fact that the absence or presence of a single edge can change $\cov_\ell(T)$ by as much as $2\ell$ prevents us from using Talagrand's inequality to show that $\cov_\ell(T)$ is tightly concentrated when $\ell$ is larger than $t^{1/3}$. Nonetheless, we believe that $\ell=o(\sqrt{t})$ is sufficient for tight concentration.

\begin{conjecture}
If $\ell=\ell(t)$ is $o(\sqrt{t})$ as $t$ approaches infinity, then $\vert \cov_\ell(T) - \EE[\cov_\ell]\vert < \delta t/\ell$ with high probability for all $\delta >0$.
\end{conjecture}

\section*{Acknowledgements}
The first author was supported by the Slovenian Research and Innovation Agency (ARIS) under the grants Z1-50003, P1-0297, and N1-0285, and by the European Union (ERC, KARST, 101071836).
The second author is funded by EPSRC (Engineering and Physical Sciences Research Council) and by the Cambridge Commonwealth, European and International Trust.
The third author was funded through a Trinity External Researcher Studentship by Trinity College of the University of Cambridge.

\printbibliography

\appendix
\section{Proof of \Cref{thm:main-queries}}
\label{ap:proof}

We first recall the following result.

\begin{lemma}[{\cite[Lemma 6]{ferber-2017}}]
    \label{lem:shorter-paths}
    Let $P=(V,E)$ be a path of length $\ell$ and let $B \subseteq E$, $|B| \leq \alpha \ell$, where $\alpha \geq \frac{1}{\ell}$. Let $Q$ be a graph obtained from $P$ by removing the edges from $B$. Then there exist vertex-disjoint subpaths $\{Q^j\}_{j \in J}$ of $Q$ such that each $Q^j$ is of length at least $\frac{1}{3 \alpha}$ and they cover at least $(\frac13 - \alpha) \ell$ vertices of $V$.
\end{lemma}

\begin{proof}[Proof of \Cref{thm:main-queries}]
    By \Cref{thm:CoverDisjointPaths}, there exist $D>1,\eps_0'>0$ such that for all $\delta\in (0,1]$ and $\eps\in (0,\eps_0')$ the following holds. With high probability, $G\sim G(n,\frac{1+\eps}{n})$ contains no system of vertex-disjoint paths of length at least $\frac{\delta}{3\eps}$ that covers more than $D\eps^2 \delta^{-2} n$ vertices.
    
    For $\delta, q\in (0,1]$,  let $\eps_0$ be the minimum of $\eps_0'$ and $\frac{q \delta^2}{192 D}$ and let $p=\frac{1+\eps}{n}$.
    Suppose that Alg is an adaptive algorithm that reveals a path of length $\ell\geq \delta/\eps$ with probability at least $q$ by querying at most $\frac{\delta q \ell}{1152 D p \eps}$ pairs of vertices in a graph $G \sim G(n,p)$. 

    Let $n' = (1 + \frac{96 D \eps^2}{q \delta^2})n$, and note that this is at most $(1+\eps/2)n$, as $\eps \leq \frac{q \delta^2}{192 D}$. Let further $V_0 = [n']$, $I_0 = \emptyset$, and $s = \frac{96 D \eps^2 n}{\delta^2 q (\ell + 1)}$. We will now construct a sequence of sets $V_0\supset V_1\supset \ldots \supset V_s$ recursively, by repeating the following steps for $i \in \{1, \ldots, s\}$.
    \begin{itemize}
        \item Sample an injection $f_i\colon [n]\rightarrow V_{i-1}$ uniformly at random and sample $G_{i} \sim G(V_{i-1}, p)$. We let $f_i\inv(G_{i})$ denote the graph with vertex set $[n]$ and edge set $\{uv\in [n]^{(2)}:f(uv)\in E(G_{i})\}$. We run Alg on $f_i\inv(G_i)$, noting that since $f_i\inv(G_i)\sim G(n,p)$, the probability of success is at least $q$.
        \item Let $L_i$ be the graph on $V_{i-1}$ with edge set $\{f(uv):uv \text{ was queried by Alg}\}$ (note that $|E(L_i)| \leq \frac{q \ell \delta}{1152 D p \eps}$), and let $K_i \subseteq L_i$ be the intersection of $L_i$ and $G_i$.
        \item If $K_i$ contains a path of length $\ell$, then we fix one such path which we call $P_i$ and let $V_i = V_{i-1} \setminus V(P_i)$ and $I_i = I_{i-1} \cup \{i\}$. Otherwise, we let $V_i = V_{i-1}$ and $I_i = I_{i-1}$.
    \end{itemize}
    Observe that we are indeed able to repeat these steps $s$ times since for every $i \in [s]$, $|V_{i-1}| \geq |V_s| \geq n' - (\ell+1) s = n$. 

    We now define a graph $H$ on the vertex set $V_0$ so that $\{u,v\} \in E(H)$ if and only if $\{u,v\} \in E(L_i)$ for some $i \in [s]$ and for the smallest such $i_0$ it holds that $\{u,v\} \in E(K_{i_0})$. 
    For every $e\in [n']^{(2)}$, we have
    \begin{equation*}
        \prob(e\in E(H))\leq \frac{1+\eps}{n} = \frac{1+\eps}{n'}+\frac{1+\eps}{n'}\left(\frac{n'}{n}-1\right)\leq \frac{1+2\eps}{n'},
    \end{equation*}
    and while the presence of edges in $H$ is not necessarily independent, one can check that for any distinct $e_1,\ldots,e_r\in [n']^{(2)}$,
    \begin{equation*}
        \prob(e_1,\ldots,e_r \in E(H)) \leq \left(\frac{1+2\eps}{n'}\right)^r.
    \end{equation*}    
     Thus, since the property of having a family of vertex-disjoint paths with length at least $\ell$ that cover $M$ vertices is monotone for all $\ell$ and $M$, the following holds. If $H$ contains a set of vertex-disjoint paths with length at least $\frac{\delta}{\eps}$ that cover at least $4 D \eps^2 \delta^{-2} n' = D (2 \eps)^2 \delta^{-2} n'$ vertices with probability at least $\frac{q^2}{4}$, the same holds with probability at least $\frac{q^2}{4}$ for $G\sim G(n', \frac{1+2\eps}{n'})$.

    Thus, it is enough to prove that $H$ has this property with probability at least $\frac{q^2}{4}$ since this will then be in contradiction with \Cref{thm:CoverDisjointPaths} for sufficiently large $n$. 

    For every $i \in I_s$ let $H_i$ be the graph on vertices $V_{i-1}$ with edges $\left( \bigcup_{j=1}^{i-1} E(L_j) \right) \cap V_{i-1}^{(2)}$.
    Observe that 
    \begin{align}
        \label{eq:H_i}
        |E(H_i)| \leq \frac{\eps}{6 \delta} \binom{|V_{i-1}|}{2}.
    \end{align}

    By definition, $V_{i-1} \setminus V_i$ is a path $P_i$ in $K_i$ for every $i \in I_s$. Thus, we can define $B_i = E(P_i) \cap E(H_i)$ and $Q_i$ to be the graph obtained from $P_i$ by deleting all edges in $B_i$. Clearly, $E(Q_i) \subseteq E(H)$, and the graphs $\{Q_i\}_{i \in I_s}$ are vertex-disjoint.

    Let $I = \{i \in I_s : |B_i| \leq \frac{\eps}{ \delta} \ell\}$. Applying \Cref{lem:shorter-paths} with $\alpha=\delta/\eps$, we get that for every $i \in I$, there exist vertex-disjoint subpaths $\{Q_i^j\}_{j\in J_i}$ of $Q_i$ such that each $Q_i^j$ has length at least $\frac{\delta}{3\eps}$ and they cover at least $(\frac13 - \frac{\eps}{ \delta}) \ell \geq \frac14 (\ell + 1)$ vertices of $P_i$, where the last inequality holds since $\eps<\frac{q \delta^2}{192 D}$, and thus, $\ell\geq 192$.

    If $|I| \geq \frac{sq}{3}$, then $\{Q_i^j\}_{j \in J_i}$ is a collection of paths of length at least $\frac{\delta}{3\eps}$ that cover at least $\frac14 (\ell+1) |I| \geq 4 D \eps^2 \delta^{-2} n'$. So it suffices to show that $|I| \geq \frac{sq}{3}$ with probability at least $\frac{q^2}{4}$. 

    Let $I' = [s] \setminus I$. For every $i \in [s]$, we have 
    \begin{align}
    \label{eq:I'}
    \prob(i \in I') & = \prob(i \notin I_s) + \prob\left(i \in I' \mid i \in I_s\right) \prob(i \in I_s).
    \end{align}
    Recall that $\prob(i \in I_s)$ is simply the success probability of Alg, which is at least $q$. For $i \in I_s$, consider the random embedding $f_i\colon [n]\rightarrow V_{i-1}$ that was chosen to pull $G_i$ back onto $[n]$. Since $f_i$ was sampled independently from all random variables that determined $H_i$ and $P_i$ is the image of a path under $f_i$, we have by~\eqref{eq:H_i}, that $\prob(e \in E(H_i)) \leq \frac{\eps}{6 \delta}$ for every $e \in E(P_i)$. By the linearity of expectation we have that $\EE[|B_i|] = \EE[|E(P_i) \cap E(H_i)|] \leq \frac{\eps}{6 \delta} \ell$. Using Markov's inequality, we obtain that 
    \begin{equation*}
        \prob(i \in I' \mid i \in I_s) = \prob\left(|B_i| \geq \frac{\eps}{\delta} \ell \mid i \in I_s\right) \leq \frac{\EE[|B_i| \mid i \in I_s]}{\frac{\eps}{\delta} \ell} \leq \frac{\frac{\eps}{6 \delta} \ell}{\frac{\eps}{\delta} \ell} < \frac{1}{2}.
    \end{equation*}
     Using these inequalities in~\eqref{eq:I'} yields $$\prob(i \in I') \leq 1 - \prob(i \in I_s) + \frac12 \prob(i \in I_s) = 1 - \frac12 \prob(i \in I_s) \leq 1 - \frac{q}{2}.$$ 

    Linearity of expectation now gives $\EE[|I'|] \leq s (1 - \frac{q}{2})$, and Markov's inequality implies that $\prob(|I'| \geq \frac{s}{1 + \frac{q}{2}}) \leq 1 - \frac{q^2}{4}$. Therefore, as $q \in (0,1)$, $$\prob(|I| \geq \frac{sq}{3}) \geq \prob(|I| \geq \frac{sq}{2+q}) = \prob(|I'| \leq s - \frac{sq}{2+q}) = \prob(|I'| \leq \frac{s}{1+\frac{q}{2}}) \geq \frac{q^2}{4},$$ as desired.
\end{proof}

\end{document}